\documentclass{psapm-l}
\usepackage{bm, stmaryrd, amscd}
\newtheorem{theorem}{Theorem}[section]
\newtheorem{lemma}[theorem]{Lemma}

\theoremstyle{definition}
\newtheorem{definition}[theorem]{Definition}

\theoremstyle{remark}

\numberwithin{equation}{section}

\newcommand{\vrho}{\varrho}
\newcommand{\Om}{\Omega}
\newcommand{\R}{\mathbb{R}}

\newcommand{\Dom}{(0,T)\times\Omega}
\newcommand{\mcw}{\mathcal{W}}
\newcommand{\Dt}{\Delta t}

\newcommand{\In}{\textrm{ in }}
\newcommand{\abs}[1]{\lvert#1\rvert}
\newcommand{\norm}[1]{\left\| #1 \right\|}
\newcommand{\Set}[1]{\left\{#1 \right\}}
\newcommand{\jump}[1]{\llbracket #1 \rrbracket}
\newcommand{\vc}[1]{\bm{#1}}
\newcommand{\weak}{\rightharpoonup}

\newcommand{\Div}{\operatorname{div}}
\newcommand{\Curl}{\operatorname{curl}}

\begin{document}

\title[]{Convergent finite element methods for compressible barotropic Stokes systems}

%    Information for authors
\author{Kenneth H. Karlsen}
\author{Trygve K. Karper}

\address{Center of Mathematics for Applications (CMA), University of Oslo, 
P.O. Box 1053, Blindern, N--0316 Oslo, Norway }
\email{kennethk@math.uio.no, t.k.karper@cma.uio.no}

% 35M10 PDE of mixed type
% 74S05 Finite element methods
% 65M12 Stability and convergence of numerical methods

\date{onsdag; 26. november, 2008}

\thanks{This work  was supported by the Research Council of Norway through 
an Outstanding Young Investigators Award.  This article was written as part 
of the the international research program 
on Nonlinear Partial Differential Equations at the Centre for 
Advanced Study at the Norwegian Academy of Science and Letters in 
Oslo during the academic year 2008--09.}

\begin{abstract}
We propose finite element methods for compressible barotropic Stokes
systems.  We state convergence results for these methods and outline 
their proofs. The principal tools of the proofs are higher integrability 
estimates for the discrete density, equations for the discrete effective 
viscous flux, and renormalized formulations of the numerical 
method for the density equation. 
\end{abstract}

\maketitle

\section{Introduction}
In this contribution we consider mixed type systems of the form
\begin{align}
		\partial_{t}\vrho + \Div(\vrho \vc{u}) &= 0, 
		\quad \In (0,T)\times \Om,     \label{eq:contequation} \\
		-\mu \Delta \vc{u} - \lambda D \Div \vc{u} + Dp(\vrho) &= \vc{f}, 
		\quad \In (0,T)\times \Om, \label{eq:momentumeq}
\end{align}
with initial data
\begin{align}\label{initial}
	\varrho|_{t=0} & = \varrho_{0},
	\quad \textrm{on $\Omega$}.
\end{align}
Here $\Omega$ is a simply connected, bounded, open, polygonal domain in $\R^N$ ($N=2,3$), 
with Lipschitz boundary $\partial \Omega$, and $T>0$ is a final time. 
The unknowns are the density $\varrho = \varrho(t,\vc{x}) \geq 0$ 
and the velocity $\vc{u} = \vc{u}(t,\vc{x}) \in \R^N$, with $t \in (0,T)$ and $\vc{x} \in \Omega$.
We denote by $\Div$ and $D$ the usual spatial divergence 
and gradient operators and by $\Delta$ the Laplace operator.

The pressure $p(\varrho)$ is governed by the equation of state 
$p(\varrho) = a\varrho^\gamma$, $a>0$ (Boyle's law). 
Typical values of $\gamma$ range from a maximum of $\frac{5}{3}$ for monoatomic gases,  
through $\frac{7}{5}$ for diatomic gases \emph{including air}, to lower values close to $1$ for 
polyatomic gases at high temperatures. We will assume that $\gamma \geq 1$.
Furthermore, the viscosity coefficients $\mu, \lambda$ 
are assumed to be constant and to
satisfy $\mu> 0, N\lambda + 2 \mu \geq 0$.

At the boundary $\partial \Om$, the system \eqref{eq:contequation}--\eqref{eq:momentumeq}
is supplemented either with the homogenous Dirichlet condition
\begin{equation}\label{eq:dirichlet}
	\vc{u} = 0, \quad \textrm{on } (0,T) \times \partial \Om,
\end{equation}
or with the Navier--slip condition
\begin{equation}\label{eq:navier}
	\vc{u} \cdot \nu =  0, \quad \Curl \vc{u} \times \nu = 0, 
	\quad \textrm{on }(0,T)\times \partial \Om.
\end{equation}

System \eqref{eq:contequation}--\eqref{eq:momentumeq} 
can be motivated in several ways. Firstly, it can be used 
as a model equation for the barotropic 
compressible Navier--Stokes equations.  This is a reasonable 
approximation for strongly viscous fluids for which convection can be neglected. 
Secondly,  in \cite[Section 5.2, Remark 5.8]{Lions}, Lions 
construct solutions to the barotropic compressible Navier--Stokes equations 
using solutions of the system \eqref{eq:contequation}--\eqref{eq:momentumeq}.  
Finally, by setting $\gamma =1$, $\vc{f} = 0$, and 
$\mu = 0$, the system \eqref{eq:contequation}--\eqref{eq:momentumeq} is 
exactly on the same form as the model derived in \cite{Weinan} for the 
dynamics of vortices in Ginzburg--Landau theories 
in superconductivity. 

Among many others, the semi--stationary system \eqref{eq:contequation}--\eqref{initial} 
has been studied by Lions in \cite[Section 8.2]{Lions} where he 
proves the existence of weak solutions and some higher regularity results. 

%:Current position - onsdag; 26. november, 2008, 22:05
The plan of this contribution is to summarize some results \cite{Karlsen1, Karlsen2,Karlsen3} 
from an ongoing project to develop convergent numerical methods
for multi-dimensional compressible viscous flow models. We 
construct numerical methods that comply with the mathematical framework 
developed for the compressible Navier--Stokes
equations by Lions \cite{Lions} and Feireisl \cite{Feireisl}. 
Over the years, several numerical methods appropriate 
for compressible viscous gas flow have been proposed.  
Except for some one-dimensional situations (cf.~Zhao and 
Hoff \cite{Zhao1, Zhao2}), it is not known, however, that these 
methods converge to a weak solution as the discretization parameters tend to zero. 
Convergence analysis for the compressible Navier--Stokes system is made 
difficult by the non--linearities in the convection and pressure terms and their
interaction. As a first step towards establishing convergence of numerical methods 
for the full system, we consider simplified systems that contain some of the difficulties but 
not all. In that respect \eqref{eq:contequation}--\eqref{eq:momentumeq} provides an example. 

The finite element methods presented here are designed to satisfy the properties needed 
to apply the weak convergence techniques used in the global existence theory for the 
compressible Navier--Stokes equations. Although the simplified system 
\eqref{eq:contequation}--\eqref{eq:momentumeq} contain additional 
structures rendering the solutions more regular than those of the full 
Navier--Stokes system,  we strive to employ techniques that can potentially be 
extended to the full system. More specifically, our finite element methods are designed 
such that Hodge decompositions of the velocity, $\vc{u} = \Curl \vc{\xi} + Dz$, can be achieved and described 
at the discrete level. This is important since then a discrete equation for the
effective viscous flux, $(\lambda + \mu)\Div \vc{u} - p(\vrho)$, can easily be extracted from the numerical scheme.  
It is the properties of this quantity that leads to strong convergence of the numerical density function;  the major
obstacle to proving convergence of a numerical method.  
Formally, multiplying the equation \eqref{eq:momentumeq} with $\vc{u}$, integrating by parts,
and using the continuity equation multiplied with $\frac{1}{\gamma-1}p'(\vrho)$ one obtains the energy relation
$$
\frac{d}{dt}\int_{\Om} \frac{p(\vrho)}{\gamma - 1}\ dx + \int_{\Om}\mu |D\vc{u}|^2 + \lambda|\Div \vc{u}|^2 \ dx 
= \int_{\Om}\vc{f}\vc{u}\ dx.
$$
A similar relation holds for our finite element methods, which reveals 
the rather weak a priori estimates that are available to us. Indeed, it is now clear 
that a major obstacle  is to obtain enough compactness on the numerical 
density $\vrho_{h}$ to conclude that $p(\vrho_{h}) \weak p(\vrho)$; of course, this 
is equivalent to $\vrho_{h} \rightarrow \vrho$ almost everywhere. 

The remaining part of this contribution is organized as follows: We collect some preliminary 
material, including the notion of weak solutions, in Section \ref{sec:prelim}. In 
Section \ref{sec:FEM} we present a finite element method 
for  the semi--stationary Stokes system in primitive variables. 
We state a convergence result for this method and comment on its proof. 
This method is fully developed and analyzed in \cite{Karlsen2}. 
In Section \ref{sec:FEM-alt}, we present and analyze an alternative 
finite element method \cite{Karlsen1}  for the same system. This method 
is, however, restricted to the case of the Navier--slip boundary condition \eqref{eq:navier}. 
Finally, we conclude this contribution by presenting a convergent 
finite element method for the Stokes approximation equations, which generalizes 
the system \eqref{eq:contequation}--\eqref{eq:momentumeq} by adding an additional 
time derivative term $\partial_t u$ to the equation for the velocity.

\section{Preliminary material}\label{sec:prelim}
Throughout the text we make frequent use of the divergence and 
curl operators and denote these by $\Div$ and $\Curl$, respectively.
In the 2D case we denote both the rotation operator taking 
scalars into vectors and the curl operator taking vectors into 
scalars by $\Curl$.  We make use of the spaces
\begin{align*}
	\vc{W}^{\Div, 2}(\Omega) & = 
	\Set{\vc{v} \in \vc{L}^2(\Omega): \Div \vc{v} \in L^2(\Omega)}, \\
	\vc{W}^{\Curl, 2}(\Omega) & = 
	\Set{\vc{v} \in \vc{L}^2(\Omega): \Curl \vc{v} \in \vc{L}^2(\Omega)},
\end{align*}
where $\nu$ denotes the unit outward pointing normal vector on $\partial \Omega$. 
If $\vc{v}\in \vc{W}^{\Div, 2}(\Omega)$ satisfies $\vc{v} \cdot \nu|_{\partial \Omega}=0$, we 
write $\vc{v}\in \vc{W}^{\Div, 2}_{0}(\Omega)$. Similarly, $\vc{v}\in \vc{W}^{\Curl, 2}_0(\Omega)$ 
means $\vc{v}\in\vc{W}^{\Div, 2}(\Omega)$ and $\vc{v} \times \nu|_{\partial \Omega} = 0$. 
In two dimensions, $\vc{w}$ is a scalar function and the 
space $\vc{W}^{\Curl,2}_{0}(\Omega)$ is to be 
understood as $W_{0}^{1,2}(\Omega)$. To define 
weak solutions, we shall use the space
$$
\mcw =\Set{\vc{v}\in  \vc{L}^2(\Omega): 
\Div \vc{v} \in L^2(\Omega), 
\Curl \vc{v} \in \vc{L}^2(\Omega), 
\vc{v} \cdot \nu|_{\partial \Omega}=0},
$$
which coincides with $ \vc{W}^{\Div, 2}_{0}(\Omega)\cap \vc{W}^{\Curl, 2}(\Omega)$. 
The space $\mcw$ is equipped with the 
norm $\norm{\vc{v}}_{\mcw}^2=\norm{\vc{v}}_{\vc{L}^2(\Omega)}^2
+\norm{\Div \vc{v}}_{\vc{L}^2(\Omega)}^2+\norm{\Curl \vc{v}}_{\vc{L}^2(\Omega)}^2$. 
It is known that $\norm{\cdot}_{\mcw}$ is equivalent to the $H^1$ norm 
on the space $\Set{v\in H^1(\Omega): \vc{v} 
\cdot \nu|_{\partial \Omega}=0}$.

Next we introduce the notion of weak solutions.

\begin{definition}[Weak solutions]\label{def:weak}
A pair $(\varrho,\vc{u})$ of functions constitutes a weak solution
of the semi-stationary compressible Stokes system \eqref{eq:contequation}--\eqref{eq:momentumeq} 
with initial data \eqref{initial} provided that:

\begin{enumerate}
	\item $(\varrho,\vc{u}) \in L^\infty(0,T;L^\gamma(\Omega))\times L^2(0,T;\mcw(\Om)),$
	\item $\partial_t\varrho + \Div (\varrho\vc{u}) = 0$ in the weak 
	sense, i.e, $\forall \phi \in C^\infty([0,T)\times\overline{\Omega})$,
	\begin{equation}\label{eq:weak-rho}
		\int_{0}^T\int_{\Omega}\varrho \left( \phi_{t} + \vc{u}D\phi\right)\ dxdt
		+ \int_{\Omega}\varrho_{0}\phi|_{t=0}\ dx = 0;
	\end{equation}
		
	\item $-\mu \Delta \vc{u} - \lambda D\Div\vc{u} + Dp(\varrho) = \vc{f}$ in 
	the weak sense, i.e, $\forall \vc{\phi} \in \vc{C}^\infty([0,T)\times\overline{\Omega})$ for which 
	$\vc{\phi} \cdot \nu = 0$ on $(0,T)\times \partial \Omega$,
	\begin{equation}\label{eq:weak-u}
		\int_{0}^T\int_{\Omega}\mu\Curl \vc{u} \Curl \vc{\phi} 
		+ \left[(\mu + \lambda)\Div \vc{u}-p (\varrho)\right]\Div \vc{\phi}\ dxdt 
		= \int_{0}^T\int_{\Omega}\vc{f}\vc{\phi}\ dxdt,
	\end{equation}
\end{enumerate}

Whenever the Dirichlet boundary condition \eqref{eq:dirichlet} is part of the problem, we 
require that $\vc{u} \times \nu = 0$ on $(0,T)\times \partial \Om$ in (1) and 
moreover that \eqref{eq:weak-u} holds for test 
functions satisfying $\phi = 0$ on $(0,T)\times \partial \Om$.
\end{definition}

\section{A non--conforming finite element method}\label{sec:FEM}
Following \cite{Karlsen2}, in this section we present a finite element method for 
the system \eqref{eq:contequation}--\eqref{eq:momentumeq} 
appropriate for both the Dirichlet boundary condition
\eqref{eq:dirichlet} and the Navier--slip boundary condition \eqref{eq:navier}.

For discretization of the velocity we will use the \emph{Crouzeix--Raviart} element space.
Consequently, the finite element method is non--conforming in the sense that the velocity approximation 
space is not a subspace of the corresponding continuous space, $\vc{W}^{1,2}(\Om)$.
Moreover, we will use a non--standard finite element formulation. 
More precisely, the formulation implicitly use the identity
\begin{equation}\label{eq:laplace}
	\int_{\Om}D\vc{u}D\vc{v}\ dx 
	= \int_{\Om}\Curl \vc{u}\Curl \vc{v} + \Div \vc{u}\Div \vc{v}\ dx, 
\end{equation}
valid for all $\vc{u} \in \mcw(\Om)$ satisfying any of 
the two boundary conditions \eqref{eq:navier} and \eqref{eq:dirichlet}.
However, as the method is non--conforming, this identity does not hold
discretely (as a sum over elements).  Still, at the discrete level, the form 
on the right-hand side of \eqref{eq:laplace} is used. In contrast to the standard 
situation where the form on the left--hand side of \eqref{eq:laplace}
is used,  this discretization does not converge unless additional terms
controlling the discontinuities of the velocity are added \cite{Brenner}:
$$
\sum_{\Gamma \in \Gamma_{h}}\frac{1}{|\Gamma|}\int_{\Gamma}
\jump{\vc{u}\cdot \nu}_{\Gamma}\jump{\vc{v}\cdot \nu}_{\Gamma} + \jump{\vc{u}\times \nu}_{\Gamma}
\jump{\vc{v}\times \nu}_{\Gamma}\ dS(x),
$$
where $\Gamma_{h}$ is the set of faces and $\jump{\cdot}_{\Gamma}$ denotes 
the jump over the edge $\Gamma$.

The advantage with this formulation is that it enables Hodge decompositions of the numerical
velocity field. By writing $\vc{u} = \Curl \vc{\xi} + Dz$ the Laplace 
operator can be split into a curl part and a divergence part plus certain jump terms. 
This is very convenient in the convergence analysis of the method. A discrete equation for the 
effective viscous flux can then be easily obtained. The reader 
is encouraged to consult \cite{Karlsen2} for the details

Given a time step $\Dt>0$, we discretize the time interval $[0,T]$ in 
terms of the points $t^m=m\Dt$, $m=0,\dots,M$, where we assume that $M\Dt=T$.  
Regarding the spatial discretization, we let $\{E_{h}\}_{h}$ be a shape regular 
family of tetrahedral meshes of $\Omega$,  where $h$ is the maximal diameter.  
It will be a standing assumption that $h$ and $\Delta t$ are 
related such that $\Delta t = c h$, for some constant $c$. 
%By shape regular we mean that there exists a constant 
%$\kappa > 0$ such that every $E \in E_{h}$ contains a ball of radius 
%$\lambda_{E} \geq \frac{h_{E}}{\kappa}$,where $h_{E}$ is the diameter of $E$. 
For each $h$, let $\Gamma_{h}$ denote the set of faces in $E_{h}$. 

We need to introduce some
additional notation for discontinuous Galerkin schemes. 
Concerning the boundary $\partial E$ of an element $E$, we write $f_{+}$ 
for the trace of the function $f$ achieved from within the element $E$ 
and $f_{-}$ for the trace of $f$ achieved from outside $E$. 
Concerning an edge $\Gamma$ that
is shared between two elements $E_{-}$ and $E_{+}$, we will write $f_{+}$ for 
the trace of $f$ achieved from within $E_{+}$ and $f_{-}$ for the trace
of $f$ achieved from within $E_{-}$. Here $E_{-}$ and $E_{+}$ are 
defined such that $\nu$ points from $E_{-}$ to $E_{+}$, where $\nu$ is 
fixed (throughout) as one of the two possible 
normal components on each edge $\Gamma$ throughout the discretization.
We also write $\jump{f}_{\Gamma}= f_{+} - f_{-}$ for the jump of $f$ 
across the edge $\Gamma$, while forward time-differencing 
of $f$ is denoted by $\jump{f^m} = f^{m+1} - f^m$ and 
$d_{t}^h[f^m] = \frac{\jump{f^m}}{\Dt}$.

We will approximate the density in the space of piecewise constants on $E_{h}$ and
we denote this space by $Q_{h}(\Om)$. For approximation of the velocity we will 
use the \emph{Crouzeix--Raviart} \cite{CR} element space
$$
\vc{V}_{h}(\Om) = \Set{\vc{v}_{h}; \vc{v}_{h}|_{E} \in \mathcal{P}_{1}^N(E), \ 
\forall E \in E_{h},\ \int_{\Gamma}\jump{\vc{v}_{h}}\ dS(x) = 0, \ 
\forall \Gamma \in \Gamma_{h}}.
$$
To incorporate boundary conditions, we let degrees of freedom of $\vc{V}_{h}(\Om)$ 
vanish at the boundary. That is, for Navier boundary condition \eqref{eq:navier} we require 
\begin{equation*}
	\int_{\Gamma}\vc{v}_{h}\cdot \nu\ dS(x) = 0, 
	\quad \forall \Gamma \in \Gamma_{h}\cap \partial \Om, \quad \forall \vc{v}_{h} \in \vc{V}_{h},
\end{equation*}
and for the Dirichlet boundary condition \eqref{eq:dirichlet}, 
\begin{equation*}
	\int_{\Gamma}\vc{v}_{h}\ dS(x) = 0, \quad \forall \Gamma \in 
	\Gamma_{h}\cap \partial \Om, \quad \forall \vc{v}_{h} \in \vc{V}_{h}.
\end{equation*}

To the space $\vc{V}_{h}(\Om)$ we associate the semi--norm
\begin{equation*}
\begin{split}
|\vc{v}_{h}|_{\vc{V}_{h}(\Om)}^2 &= \|\Curl_{h} \vc{v}_{h}\|_{L^2(\Om)}^2 + \|\Div_{h} \vc{v}_{h}\|_{L^2(\Om)}^2 \\
&\qquad + 
\frac{h^\epsilon}{|\Gamma|}\sum_{\Gamma \in \Gamma_{h}}
\|\jump{\vc{v}_{h} \cdot \nu}\|_{L^2(\Gamma)}^2 
+ \|\jump{\vc{v}_{h} \times \nu}\|_{L^2(\Gamma)}^2,
\end{split}
\end{equation*}
and the corresponding norm
\begin{equation*}
\|\vc{v}_{h}\|^2_{\vc{V}_{h}(\Om)} = \|\vc{v}_{h}\|_{L^2(\Om)}^2 + \abs{\vc{v}_{h}}_{\vc{V}_{h}(\Om)}^2.
\end{equation*}
Here, $\Curl_{h}$ and $\Div_{h}$ denotes the curl and divergence operators, respectively,  taken inside each element.
The scaling parameter $\epsilon>0$ is required to prove
convergence of the finite element method. The size of $\epsilon$ will affect the accuracy of 
the method and it should therefore be fixed very small in practical computations  \cite{Karlsen2}.

Before stating the finite element method, we recall
from \cite{Karlsen2} the following basic compactness result for
approximations in $\vc{V}_{h}(\Om)$.
\begin{lemma}\label{lemma:compactembedding}
There exists a constant $C>0$, depending only on the shape regularity of $E_{h}$ and the size of $\Omega$,
such that for any $\xi \in \mathbb{R}^2$
\begin{equation*}%\label{eq:compact}
	\|\vc{v}_{h}(\cdot) - \vc{v}_{h}(\cdot-\xi)\|_{\vc{L}^2(\Om)} 
	\leq C|\xi|^{\frac{1}{2}- \frac{\epsilon}{4}}|\vc{v}_{h}|_{\vc{V}_{h}(\Om)}, 
	\quad \forall \vc{v}_{h} \in \vc{V}_{h}(\Om),
\end{equation*}
and 
%\begin{equation*}%\label{eq:poincare}
	$\|\vc{v}_{h}\|_{\vc{L}^2(\Om)} \leq C|\vc{v}_{h}|_{\vc{V}_{h}(\Om)}$, 
	%\quad 
	$\forall \vc{v}_{h} \in \vc{V}_{h}(\Om)$.
%\end{equation*}
\end{lemma}

%The finite element method is given in the following definition.

\begin{definition}[Finite element method]\label{def:num-schemea}
Let $\Set{\vrho^0_h(x)}_{h>0}$ be a sequence (of piecewise constant 
functions) in $Q_{h}(\Omega)$ that satisfies $\vrho_h^0>0$ for each fixed $h>0$ 
and $\varrho^0_h\to \vrho^0$ a.e.~in $\Om$ and in $L^1(\Om)$ 
as $h\to 0$. Set $\vc{f}_{h}:=\Pi_{h}^{Q} \vc{f}$, where it is understood that 
$\Pi_{h}^{Q} \vc{f}$ projects $\vc{f}(t,x)$ onto constants both in time $t$ and space $x$; for notational 
convenience we set $\vc{f}_h^m:=\vc{f}_h(t^m,\cdot)\in Q_{h}(\Omega)$ 
for any $m=0,\dots,M$.

Now, determine functions 
$(\varrho^m_{h},\vc{u}^m_{h}) \in Q_{h}(\Omega)
\times \vc{V}_{h}(\Omega)$, $m=1,\dots,M$, such that 
for all $\phi_{h} \in Q_{h}(\Omega)$,
\begin{equation}\label{FEMp:contequation}
	\begin{split}
		&\int_\Omega d_{t}^h[\varrho^m_h] \phi_{h}\ dx
		-\Delta t\sum_{\Gamma \in \Gamma_h}\int_\Gamma \left(\varrho^m_{-}(\vc{u}^{m}_h \cdot \nu)^+
		+\varrho^m_+(\vc{u}^{m}_h \cdot \nu)^-\right)\jump{\phi_{h}}_\Gamma\ dS(x)
		%\\ & \qquad 
		= 0.
	\end{split}
\end{equation}
and for all $\vc{v}_{h} \in  \vc{V}_{h}(\Omega)$,
\begin{equation}\label{FEMp:momentumeq}
	\begin{split}
		&\int_{\Om} \mu\Curl_{h} \vc{u}^m_{h}\Curl_{h} \vc{v}_{h} + \left[(\mu + \lambda)\Div_{h} \vc{u}^m_{h}
		-p(\varrho^m_{h})\right]\Div_{h} \vc{v}_{h}\ dx \\
		&\quad +\sum_{\Gamma \in \Gamma_{h}}\frac{h^\epsilon}{|\Gamma|}\int_{\Gamma}
		\jump{\vc{u}^m_{h}\cdot \nu}\jump{\vc{v}_{h}\cdot \nu}
		+ \jump{\vc{u}^m_{h}\times \nu}\jump{\vc{v}_{h}\times \nu}\ dS(x) 
		%\\ & \qquad 
		= \int_{\Omega}\vc{f}^m_{h}\vc{v}_{h}\ dx, 
	\end{split}
\end{equation}
%for $m=1,\dots,M$. 

In \eqref{FEMp:contequation}, $(\vc{u}_{h} \cdot \nu)^+(x)
=\max\left\{\frac{1}{|\Gamma|}\int_{\Gamma} \vc{u}_{h}\cdot \nu\ dS(x) , 0\right\}$ 
and \\ $(\vc{u}_{h} \cdot \nu)^-(x) = \min\left\{\frac{1}{|\Gamma|}\int_{\Gamma} \vc{u}_{h} 
\cdot \nu\ dS(x), 0\right\}$ for $x \in \Gamma$ 
and all $\Gamma \in \Gamma_{h}$.
\end{definition}

The existence of a solution to the discrete 
equations \eqref{FEMp:contequation}--\eqref{FEMp:momentumeq} is proved 
in \cite{Karlsen2} by using a topological degree argument.
In \cite{Karlsen2} it is also shown that the scheme preserves the total mass 
and that the density remains strictly positive provided that the initial density is strictly positive.
Moreover, for any $m=1, \ldots, M$, 
\begin{equation*}%\label{eq:stabilityincerrors}
	\begin{split}
		&\int_{\Omega} P(\varrho_{h}^m)\ dx 
		+ C \sum_{k=1}^m \Delta t \|\vc{u}_{h}^k\|_{\vc{V}_{h}(\Omega)}^2 \\
		& \quad 
		+ \sum_{k=1}^{m}\int_{\Omega} P''(\varrho_{\dagger\dagger}^k)\jump{\varrho^{k-1}_{h}}^2\ dx
		+ \sum_{k=1}^{m}\sum_{\Gamma \in \Gamma_{h}}\Delta t
		\int_{\Gamma}P''(\varrho^k_{\dagger})\jump{\varrho^k_{h}}_{\Gamma}^2
		\abs{\vc{u}^k_{h} \cdot \nu}\ dx \\
		& \leq \int_{\Omega} P(\varrho_{0}) \ dx 
		+\frac{1}{4C} \sum_{k=1}^m\Delta t\|\vc{f}_{h}^k\|_{\vc{L}^2(\Om)}^2,
	\end{split}
\end{equation*}
where $P(\varrho) = \frac{p(\varrho)}{\gamma - 1}$ if $\gamma > 1$ 
and $P(\vrho) = \vrho \log \vrho$ if $\gamma =1$. 
Moreover, $\vrho_{\dagger \dagger}^k \in [\vrho^{k-1}_{h}, \vrho_{h}^k]$  
and $\vrho_{\dagger}^k \in [\vrho_{+}^k, \vrho^k_{-}]$.

Next, for each fixed $h>0$, we extend the numerical solution 
$\Set{(\varrho^m_{h},\vc{u}^m_{h})}_{m=0}^M$
to the whole of $(0,T)\times \Omega$ by setting 
\begin{equation}\label{eq:num-scheme-IIa}
	(\varrho_{h},\vc{u}_{h})(t)=(\varrho^m_{h},\vc{u}^m_{h}), 
	\qquad t\in (t_{m-1},t_m), \quad m=1,\dots,M.
\end{equation}
In addition, we set $\varrho_{h}(0)= \varrho^0_{h}$. 

The main result of \cite{Karlsen2} is that the approximate solutions \eqref{eq:num-scheme-IIa}
converge  to a weak solution of the semi--stationary 
Stokes system \eqref{eq:contequation}--\eqref{eq:momentumeq}.

\begin{theorem}\label{theorem:mainconvergence2}
Suppose $\vc{f}\in \vc{L}^2(\Dom)$ and $\vrho_0 \in L^\gamma (\Om)$, if $\gamma > 1$,
 and $\varrho_{0}\log \varrho_{0} \in L^1(\Om)$, if $\gamma=1$.  
Let $\Set{(\varrho_{h},\vc{u}_{h})}_{h>0}$ be a sequence of numerical solutions 
constructed according to \eqref{eq:num-scheme-IIa} and Definition \ref {def:num-schemea}. 
Then, passing if necessary to a subsequence as $h\to 0$, 
$\vc{u}_{h} \weak \vc{u}$ in $L^2(0,T;\vc{L}^{2}(\Omega))$, 
$\varrho_{h}\vc{u}_{h} \weak \varrho\vc{u}$ in the sense 
of distributions on $\Dom$, and $\varrho_{h} \rightarrow \varrho$
a.e.~in $\Dom$, where the limit pair $(\vrho, \vc{u})$  is a weak 
solution as stated in Definition \ref{def:weak}.
\end{theorem}

\subsubsection*{Comments on the proof of Theorem \ref{theorem:mainconvergence2}}
In proving convergence to a weak solution of the continuity 
equation \eqref{eq:contequation} the main step is to
obtain convergence of the product $\vrho_{h}\vc{u}_{h}Ê\weak \vrho \vc{u}$
in the sense of distributions; this follows from an Aubin--Lions argument
using the spatial compactness of the velocity established in Lemma \ref{lemma:spacetranslation} 
combined with the fact that $d_{t}^h[\vrho_{h}] \in_{\mathrm{b}} L^1(0,T;W^{-1,1}(\Om))$. 

To conclude convergence to a weak solution 
of the velocity equation \eqref{eq:momentumeq}, we need a 
higher integrability estimate for the numerical density. 
We achieve this by utilizing test functions $\vc{v}_{h} \in \vc{V}_{h}(\Om)$ 
satisfying $\Div \vc{v}_{h} = p(\vrho_{h})$, thereby obtaining 
$p(\vrho_{h}) \in_{\mathrm{b}} L^2(0,T;L^2(\Om))$. 
Next, we establish 
strong convergence of the density. This is obtained by first proving 
weak sequential continuity of the effective viscous flux. 
That is, first we establish that $\lim_{h \rightarrow 0}\left[(\lambda + \mu)
\Div \vc{u}_{h} - p(\vrho_{h})\right]\vrho_{h} 
= \overline{(\mu+ \lambda)\Div \vc{u} - p(\vrho)}\vrho$, where the overbar denotes 
the weak limit. In this step, the div--curl structure of the scheme is utilized. In particular, we 
employ test functions $\vc{v}_{h} \in \vc{V}_{h}(\Om)$
that satisfies $\Div \vc{v}_{h} = \vrho_{h}$ and $\Curl \vc{v}_{h} = 0$ on elements away from the boundary.
Finally, using this and a renormalized version of the continuity scheme \eqref{FEMp:contequation}, 
we obtain strong convergence of the density.

\section{A mixed finite element method}\label{sec:FEM-alt}
Following \cite{Karlsen1}, we present an alternative finite 
element method appropriate for the Navier--slip boundary 
condition \eqref{eq:navier}. The method is derived by introducing 
the vorticity $\vc{w} = \Curl \vc{u}$ as an auxiliary 
variable and recasting \eqref{eq:momentumeq} as
$$
\mu \Curl \vc{w} - (\lambda + \mu)D\Div \vc{u} + Dp(\vrho) = \vc{f},
$$
where also the identity $-\Delta = \Curl \Curl -D \Div$ is used.
This leads naturally to the following mixed formulation: Determine functions 
$$
(\vc{w},\vc{u}) \in L^2(0,T;\vc{W}_{0}^{\Curl, 2}(\Omega))
\times L^2(0,T;\vc{W}^{\Div, 2}_{0}(\Omega))
$$
such that
\begin{equation}\label{def:mixed-weak}
	\begin{split}
		& \int_{0}^T\int_{\Omega}\mu\Curl \vc{w} \vc{v}
		+ \left[(\mu + \lambda)\Div \vc{u} - p(\varrho)\right]\Div \vc{v}\ dxdt
		= \int_{0}^T\int_{\Omega}\vc{f}\vc{v}\ dxdt, \\
		& \int_{0}^T\int_{\Omega} \vc{w}\vc{\eta}-\Curl \vc{\eta} \vc{u}\ dxdt = 0, 
	\end{split}
\end{equation}
for all $(\eta,\vc{v}) \in  L^2(0,T;\vc{W}_{0}^{\Curl, 2}(\Omega))
\times L^2(0,T;\vc{W}^{\Div, 2}_{0}(\Omega))$. 
We make clear that if $(\vrho, \vc{w}, \vc{u})$ is a triple 
satisfying \eqref{eq:weak-rho} and \eqref{def:mixed-weak},
then the pair $(\vrho,\vc{u})$ is also a weak solution 
according to Definition \ref{def:weak}. 

To obtain a stable numerical method, the mixed 
finite element formulation of \eqref{def:mixed-weak} is posed
with the velocity $\vc{v}_{h}$ in a div--conforming space 
$\vc{V}_{h}(\Om) \subset \vc{W}^{\Div,2}_{0}(\Om)$ 
and vorticity $\vc{w}_{h}$ in a curl--conforming space 
$\vc{W}_{h}(\Om) \subset \vc{W}^{\Curl,2}_{0}(\Om)$.
There exists several such spaces, however here we will use 
the Nedelec spaces of first order and
first kind. We choose these spaces for their simplicity and since the 
most natural choice of approximation space for the density 
is then the space of piecewise constants.
We will continue to denote this space by $Q_{h}(\Om)$. 

This choice of finite element spaces is also very convenient since they can be 
related through the exact de Rham sequence
\begin{equation*}
	\begin{CD}
		0 @> \subset >> S_{h} @> 
		\operatorname{grad} >> \vc{W}_h @> \Curl\ >> \vc{V}_h @> 
		\Div\ >> Q_{h} @> >> 0.
	\end{CD}
\end{equation*} 
Thus, we can use spaces orthogonal to the 
range of the previous operator, i.e., 
\begin{equation*}
	\vc{W}_{h}^{0,\perp}  := \{\vc{w}_{h} \in \vc{W}_{h}; 
	\Curl \vc{w}_{h} =0 \}^\perp\cap \vc{W}_{h}, \quad 
	\vc{V}_{h}^{0,\perp}  := \{\vc{v}_{h} \in \vc{V}_{h}; 
	\Div \vc{v}_{h} = 0\}^\perp \cap \vc{V}_{h},
\end{equation*}
to deduce the decompositions
\begin{equation*}
	\vc{W}_{h} = DS_{h} +\vc{W}_{h}^{0, \perp}, \quad 
	\vc{V}_{h}  = \Curl \vc{W}_{h} +\vc{V}_{h}^{0, \perp},
\end{equation*}
together with the discrete Poincar\'e inequalities
\begin{equation*}
	\norm{\vc{v}_{h}}_{\vc{L}^2(\Omega)}  \leq 
	C\norm{\Div \vc{v}_{h}}_{L^2(\Omega)},
	\quad
	\norm{\vc{w}_{h}}_{\vc{L}^2(\Omega)} 
	\leq C\norm{\Curl \vc{w}_{h}}_{L^2(\Omega)},
\end{equation*}
Consequently, as with the previous method, the mixed 
finite element method also admits Hodge decompositions, which in turn 
implies that a discrete equation for the effective viscous flux can be derived.

We need the following compactness property of the space 
$\vc{V}^{0,\perp}_{h}$. The proof is given in \cite[Appendix A]{Karlsen1}.

\begin{lemma}\label{lemma:spacetranslation}
Let $\{\vc{v}_{h}\}_{h>0}$ be a sequence in $\vc{V}^{0,\perp}_{h}$ 
such that $\|\Div \vc{v}_{h}\|_{L^2(\Omega)}\leq C$, where 
the constant $C>0$ is independent of $h$. 
Then, for any $\xi \in \mathbb{R}^N$,
\begin{equation*}
	\|\vc{v}_{h}(x) - \vc{v}_{h}(x- \xi) \|_{\vc{L}^2(\Omega)} 
	\leq C(|\xi|^{\frac{4-N}{2}} + |\xi|^2)^\frac{1}{2}
	\|\Div \vc{v}_{h}\|_{L^2(\Omega)},
\end{equation*}
where the constant $C>0$ is independent of both $h$ and $\xi$.
\end{lemma}

\begin{definition}[Mixed finite element method]\label{def:num-scheme2}
Let $\vrho^0_h$ and  $\vc{f}_h^m$ be as given in Definition \ref{def:num-schemea}. 
Determine functions 
$$
(\varrho^m_{h},\vc{w}^m_{h},\vc{u}^m_{h}) \in Q_{h}(\Omega)\times\vc{W}_{h}(\Omega)
\times \vc{V}_{h}(\Omega), \quad m=1,\dots,M,
$$
such that for all $\phi_{h} \in Q_{h}(\Omega)$,
\begin{equation}\label{FEM:contequation2}
	\begin{split}
		&\int_\Omega d_{t}^h[\varrho^m_h] \phi_{h}\ dx
		-\Delta t\sum_{\Gamma \in \Gamma_h}\int_\Gamma \left(\varrho^m_{-}(\vc{u}^{m}_h \cdot \nu)^+
		+\varrho^m_+(\vc{u}^{m}_h \cdot \nu)^-\right)[\phi_{h}]_\Gamma\ dS(x)=0,
	\end{split}
\end{equation}
and for all $(\vc{\eta}_{h},\vc{v}_{h}) \in \vc{W}_{h}(\Omega)\times \vc{V}_{h}(\Omega)$,
\begin{equation}\label{FEM:momentumeq2}
	\begin{split}
		&\int_{\Omega}\mu\Curl \vc{w}^m_{h}\vc{v}_{h} + \left[(\mu + \lambda)\Div \vc{u}^m_{h}
		-p(\varrho^m_{h})\right]\Div \vc{v}_{h}\ dx
		= \int_{\Omega}\vc{f}^m_{h}\vc{v}_{h}\ dx, \\
		&\int_{\Omega}\vc{w}^m_{h}\vc{\eta}_{h} -  \vc{u}^m_{h}\Curl \vc{\eta}_{h} \ dx =0.
	\end{split}
\end{equation}
In \eqref{FEM:contequation2}, $(\vc{u}_{h} \cdot \nu)^+=\max\{\vc{u}_{h} \cdot \nu, 0\}$ 
and $(\vc{u}_{h} \cdot \nu)^- = \min\{\vc{u}_{h} \cdot \nu, 0\}$.
\end{definition}

The existence of a solution to the discrete 
equations \eqref{FEM:contequation2}--\eqref{FEM:momentumeq2} 
is proved in \cite{Karlsen1}. Moreover, for any $m=1, \ldots, M$,  
\begin{equation*}%\label{eq:stabilityincerrors}
	\begin{split}
		&\int_{\Omega} P(\varrho_{h}^m)\ dx 
		+\sum_{k=1}^m\Delta t \|\vc{u}_{h}^k\|^2_{\vc{W}^{\Div, 2}(\Om)}
		+\sum_{k=1}^m\Delta t\|\vc{w}_{h}^k\|^2_{\vc{W}^{\Curl,2}(\Om)} \\
		& \quad 
		+ \sum_{k=1}^{m}\int_{\Omega} P''(\vrho_{\dagger \dagger}^k)\jump{\varrho^{k-1}_{h}}^2\ dx
		+ \sum_{k=1}^{m}\sum_{\Gamma \in \Gamma_{h}}\Delta t
		\int_{\Gamma}P''(\varrho^k_{\dagger})\jump{\varrho^k_{h}}_{\Gamma}^2
		\abs{\vc{u}^k_{h} \cdot \nu}\ dx \\
		& \leq \int_{\Omega} P(\varrho_{0}) \ dx 
		+C \sum_{k=1}^m\Delta t\|\vc{f}_{h}^k\|_{\vc{L}^2(\Om)}^2,
	\end{split}
\end{equation*}
where $\vrho_{\dagger \dagger}^k \in [\vrho^{k-1}_{h}, \vrho_{h}^k]$ 
and $\vrho_{\dagger}^k \in [\vrho_{+}^k, \vrho^k_{-}]$.

For each fixed $h>0$, the numerical solution 
$\Set{(\varrho^m_{h},\vc{w}^m_{h},\vc{u}^m_{h})}_{m=0}^M$
is extended to the whole of $(0,T)\times \Omega$ by setting 
\begin{equation}\label{eq:num-scheme-IIb}
	(\varrho_{h},\vc{w}_{h},\vc{u}_{h})(t)=(\varrho^m_{h},\vc{w}^m_{h},\vc{u}^m_{h}), 
	\qquad t\in (t_{m-1},t_m), \quad m=1,\dots,M.
\end{equation}
In addition, we set $\varrho_{h}(0)= \varrho^0_{h}$. The main 
result in \cite{Karlsen1} is that the sequence 
$\{\varrho_{h}, \vc{w}_{h}, \vc{u}_{h}\}_{h>0}$ converges
to a weak solution in the sense of Definition \ref{def:weak}.

\begin{theorem}\label{theorem:mainconvergence1}
Suppose $\vc{f}\in \vc{L}^2(\Dom)$, $\vrho_0\in L^\gamma(\Om)$ if $\gamma > 1$, 
and $\varrho_{0}\log \varrho_{0} \in L^1(\Omega)$ if $\gamma=1$.  
Let $\Set{(\varrho_{h},\vc{w}_{h},\vc{u}_{h})}_{h>0}$ be a sequence of numerical solutions 
constructed according to \eqref{eq:num-scheme-IIb} and Definition \ref {def:num-scheme2}. 
Then, passing if necessary to a subsequence as $h\to 0$, $\vc{w}_{h} \weak \vc{w}$ 
in $L^2(0,T;\vc{W}^{\Curl,2}_{0}(\Om))$, $\vc{u}_{h} \weak \vc{u}$ 
in $L^2(0,T;\vc{W}^{\Div, 2}_{0}(\Omega))$, 
$\varrho_{h}\vc{u}_{h} \weak \varrho\vc{u}$ in the sense 
of distributions on $\Dom$, and $\varrho_{h} \rightarrow \varrho$
a.e.~in $\Dom$, where the limit triplet $(\vrho,\vc{w}, \vc{u})$ satisfies the mixed form
\eqref{def:mixed-weak}, and consequently $(\vrho,\vc{u})$ is a weak 
solution as stated in Definition \ref{def:weak}.
\end{theorem}

\subsubsection*{Comments on the proof of Theorem \ref{theorem:mainconvergence1}}

The proof of convergence to a weak solution of the continuity equation \eqref{eq:contequation}
is similar to the corresponding step in the proof of Theorem \ref{theorem:mainconvergence2}.  
The difference is that the compactness of the velocity approximation 
now requires a different argument.  In particular, Lemma \ref{lemma:spacetranslation} must 
be employed. The proof of convergence to a weak solution 
of the velocity equation \eqref{eq:momentumeq} is also similar 
to the proof of Theorem \ref{theorem:mainconvergence2}.
However, a difference is that the weak sequential continuity of 
the effective viscous flux can now be 
obtained by using test functions in $\vc{V}_{h}^{0,\perp}(\Om)$ 
satisfying $\Div \vc{v}_{h} = \vrho_{h}$.
Strong convergence of the density is then obtained 
as in the proof of Theorem \ref{theorem:mainconvergence2}.

\section{Extension to the Stokes approximation equations.}
In this final section we present an extension of the previous finite element method to the 
following system
\begin{align}
	\partial_{t}\varrho + \Div (\varrho \vc{u}) &=0,\quad \In \Dom \label{eq:stoke1} \\
	\overline{\varrho}\partial_{t}\vc{u} - \mu \Delta \vc{u} - \lambda D \Div \vc{u} 
	+ Dp(\varrho) &= 0, \quad \In \Dom,
	\label{eq:stoke2}
\end{align}
where $\overline{\varrho} = \frac{1}{|\Om|}\int_{\Om}\varrho_{0}\ dx$  denotes the average initial density. 
The equations \eqref{eq:stoke1}--\eqref{eq:stoke2} is known in the 
literature as the \emph{Stokes approximation equations}.
The system is almost identical to the compressible Stokes 
system \eqref{eq:contequation}--\eqref{eq:momentumeq}, the
difference being the inclusion of the time derivative term in \eqref{eq:stoke2}.

The finite element method is similar to the mixed method of Section \ref{sec:FEM-alt}, and as such 
it is only applicable to the case of Navier--slip boundary conditions \eqref{eq:navier}.
Furthermore, for technical reasons, convergence is proved only in the case $\gamma > \frac{N}{2}$.
The method is constructed and analyzed in \cite{Karlsen3}.

\begin{definition}[Numerical scheme]\label{def:num-scheme}
Let $\vrho^0_h$ be as given in Definition \ref{def:num-schemea}. 
Determine functions 
$$
(\varrho^m_{h},\vc{w}^m_{h},\vc{u}^m_{h}) \in Q_{h}(\Omega)\times\vc{W}_{h}(\Omega)
\times \vc{V}_{h}(\Omega), \quad m=1,\dots,M,
$$
such that for all $\phi_{h} \in Q_{h}(\Omega)$,
\begin{equation}\label{FEM:contequation3}
	\begin{split}
		&\int_\Omega d_{t}^h[\varrho^m_h] \phi_{h}\ dx
		-\Delta t\sum_{\Gamma \in \Gamma_h}\int_\Gamma \left(\varrho^m_{-}(\vc{u}^{m}_h \cdot \nu)^+
		+\varrho^m_+(\vc{u}^{m}_h \cdot \nu)^-\right)\jump{\phi_{h}}_\Gamma\ dS(x) = 0,
	\end{split}
\end{equation}
and for all $(\vc{\eta}_{h},\vc{v}_{h}) \in \vc{W}_{h}(\Omega)\times \vc{V}_{h}(\Omega)$,
\begin{equation}\label{FEM:momentumeq3}
	\begin{split}
		&\int_{\Omega}d_{t}^h[\vc{u}_{h}^m]\vc{v}_{h} + \mu\Curl \vc{w}^m_{h}\vc{v}_{h} + \left[(\mu + \lambda)\Div \vc{u}^m_{h}
		-p(\varrho^m_{h})\right]\Div \vc{v}_{h}\ dx
		= 0, \\
		&\int_{\Omega}\vc{w}^m_{h}\vc{\eta}_{h} -  \vc{u}^m_{h}\Curl \vc{\eta}_{h} \ dx =0,
	\end{split}
\end{equation}
for $m=1,\dots,M$. 
\end{definition}

Existence of a numerical solution and various properties of these solutions hold as in the 
previous section with only minor modifications.  
We extend $\{\vrho_{h}^k, \vc{w}_{h}^k, \vc{u}_{h}^k\}_{h>0}$ for $k=1, \ldots, M$ to 
functions $\{(\varrho_{h},\vc{w}_{h},\vc{u}_{h})\}_{h> 0}$ defined on all of $(0,T)\times \Om$
as in \eqref{eq:num-scheme-IIb}. The main result in \cite{Karlsen3} is that the 
sequence $\{(\varrho_{h},\vc{w}_{h},\vc{u}_{h})\}_{h>0}$ converges to a weak solution 
of the Stokes approximation equations \eqref{eq:stoke1}--\eqref{eq:stoke2}.
The notion of a weak solution is similar to that in Definition \ref{def:weak}.

\begin{theorem}\label{theorem:mainconvergence}
Suppose  $\gamma > \frac{N}{2}$ and $\vrho_0\in L^\gamma(\Om)$ .
Let $\Set{(\varrho_{h},\vc{w}_{h},\vc{u}_{h})}_{h>0}$ be a sequence of numerical solutions 
constructed according to Definition \ref {def:num-scheme}. 
Then, passing if necessary to a subsequence as $h\to 0$, $\vc{w}_{h} \weak \vc{w}$ 
in $L^2(0,T;\vc{W}^{\Curl,2}_{0}(\Om))$, $\vc{u}_{h} \weak \vc{u}$ 
in $L^2(0,T;\vc{W}^{\Div, 2}_{0}(\Omega))$, 
$\varrho_{h}\vc{u}_{h} \weak \varrho\vc{u}$ in the sense 
of distributions, and $\varrho_{h} \rightarrow \varrho$
a.e.~in $\Dom$, where the limit triplet $(\vrho,\vc{w}, \vc{u})$ is a weak 
solution to the Stokes approximation equations \eqref{eq:stoke1}--\eqref{eq:stoke2}.
\end{theorem}

\subsubsection*{Comments to the proof of Theorem \ref{theorem:mainconvergence}}
The proof of convergence is similar to the proof of Theorem \ref{theorem:mainconvergence2}.
However, the proof of higher integrability on the density and the proof of
weak sequential continuity of the effective viscous flux now requires additional arguments in order 
to handle the time derivative in \eqref{FEM:momentumeq3}. In particular, the continuity scheme \eqref{FEM:contequation3} needs to be used to handle the term $\vc{u}\vc{v}_{t}$, where 
$\vc{v}_{t} \in \vc{V}^{0,\perp}(\Om)$ satisfies $\Div \vc{v}_{t}  = \vrho_{t}$.
For technical reasons, we must then require $\gamma> \frac{N}{2}$. Moreover,
the higher integrability estimate for the density now 
gives $p(\vrho_{h})\vrho_{h} \in_{\mathrm{b}} L^1(0,T;L^1(\Om))$.
In addition,  in this case we need in fact strong convergence of the velocity, $\vc{u}_{h} \rightarrow \vc{u}$.
This is obtained through an Aubin--Lions argument using the spatial compactness on the 
velocity together with weak control of $d_{t}^h[\vc{u}_{h}]$. 
Due to space limitations we refer the reader to \cite{Karlsen3} for details.

\bibliographystyle{amsalpha}

\end{document}